\documentclass[11pt]{amsart}

\usepackage[parfill]{parskip}    
\usepackage{graphicx, txfonts}
\usepackage{epstopdf}
\usepackage{hyperref}
\usepackage{amsthm, amssymb, amsfonts, amsmath, enumerate, xy, mathrsfs}
\input xy
\xyoption{all}
\newtheorem{theorem}{Theorem}[section]

\theoremstyle{definition}

\newtheorem{conjecture}[theorem]{Conjecture}

\newtheorem{example}[theorem]{Example}

\newtheorem{keyidea}[theorem]{Key Idea}
\newtheorem{organizing}[theorem]{Organizing Principle}

\theoremstyle{remark}

\newcommand{\N}{\mathbb{N}}

\newcommand{\M}{\mathcal{M}}

\DeclareMathOperator{\Sym}{Sym}

\newcommand{\po}{\ar@{}[dr]|(.7){\Searrow}}
\newcommand{\pb}{\ar@{}[dr]|(.3){\Nwarrow}}

\newcommand{\boxprod}{\mathbin\square}

\title{Monoidal Bousfield Localizations and Algebras over Operads: A User's Guide}

\author{David White}
\address{Department of Mathematics \\ Denison University
\\ Granville, OH 43023}
\email{david.white@denison.edu}

\date{September 1, 2015}

\begin{document}

\begin{abstract}
This paper is a companion for the paper ``Monoidal Bousfield Localizations and Algebras over Operads'' \cite{white-localization}, part of the author's PhD thesis. This paper was written in 2015 for the first edition of {\em Enchiridion: Mathematics User's Guides}, published under the banner of the {\em Journal of Humanistic Mathematics} \cite{user-guide-project}. More User's Guides can be found at \href{https://mathusersguides.com/}{https://mathusersguides.com/}.
\end{abstract}

\maketitle

\tableofcontents

\section{Key insights and central organizing principles}

A key task of mathematics is to provide unifying perspectives. Such perspectives help us build bridges between seemingly disparate fields, work in a general setting so that our results hold in numerous specific instances, and understand the humanistic and aesthetic side of our work, i.e. determine what drives us to research the topics we choose, which deep results we are truly pursuing, and which ideas the human mind keeps returning to in different guises and different times. One such idea is that of localization, which allows us to zoom in on the pertinent information within a problem. This is like putting on glasses to change which objects we view as equivalent. Another fundamental idea is 

\begin{organizing} \label{keyidea:replace}
Whenever possible, one should work in a setting where it is possible to replace the objects of interest by nicer objects which are equivalent in a suitable sense.
\end{organizing}

My paper \cite{white-localization} is fundamentally about localization. The purpose is to understand how the types of localization which arise in homotopy theory (specifically, Bousfield localization of model categories) interacts with algebraic structure in these categories. The \textbf{main goal} of the paper is to find checkable conditions so that Bousfield localization preserves algebras over operads (i.e. vehicles for encoding algebraic structure in general monoidal categories, such as commutative structure, associative structure, Lie structure, etc). En route we provide conditions so that the model category obtained via Bousfield localization satisfies various axioms that are common in monoidal model category theory. This includes conditions so that commutative monoids inherit a model structure and are preserved under localization. These conditions are then checked in numerous model categories of interest, and general results are translated into results specific to these settings, recovering several classical results (e.g. about spaces and chain complexes) in a unified framework and proving new results about equivariant spectra, ideals of ring spectra, and different models for monoidal stable homotopy theory.

In this paper localization arises in two different but related ways. In the first, localization is applied to a model category in order to obtain its homotopy category (a good example to keep in mind is the mental shift you do when thinking of topological spaces up to homeomorphism vs. up to homotopy equivalence, i.e. up to continuous deformation). This form of localization goes back to \cite{gabriel-zisman} at least, and can be viewed as a generalization of the localization which arises in algebra. There, one localizes by formally adjoining multiplicative inverses to a specified set of elements in a ring.  In the category theoretic version one does not have elements to invert, so instead one formally adjoins morphisms which are inverses to a specified class of maps. To recover the ring-theoretic notion of localization one inverts endomorphisms of the ring which correspond to multiplying by the elements one seeks to invert.

Unfortunately, not all choices of sets of maps admit localization. The notion of a model category (which comes equipped with a chosen class of maps to invert called the weak equivalences) arose to fix this issue and to provide control over the morphisms in the localized category. Having the structure of a model category allows for the tools of homotopy theory to be applied, and in this way parts of homological algebra, algebraic geometry, representation theory, logic, graph theory, and even computer science can be viewed as special cases of homotopy theory. To summarize

\begin{organizing} \label{keyidea:build-model-categories}
In settings where one has a notion of weak equivalence or something like a homology theory to compress complicated information into simple information, one should try to build a model structure so that the tools of abstract homotopy theory can be applied.
\end{organizing}

The other type of localization in this paper is called Bousfield localization, and it is a procedure one applies to a model category $\M$ in order to enlarge the specified class of weak equivalences to contain some specified set of maps $C$ (the resulting model category is denoted $L_C(\M)$), see \cite{bousfield-localization-spaces-wrt-homology}, \cite{bous79}, and \cite{hirschhorn}. Both forms of localization can be viewed as special cases of Organizing Principle \ref{keyidea:replace}, because both satisfy a universal property saying they are the ``closest'' to the given category in which the prescribed maps have been inverted. Formally, this means any functor out of the given category which inverts the maps factors through the localization. The notions of monoidal model categories, operads, and Bousfield localization are recalled in Section 2. 

This paper began out of a desire to understand an example of a localization which destroys equivariant commutativity (Example 5.7). This example arose during the recent proof of Hill, Hopkins, Ravenel \cite{kervaire-arxiv} of the Kervaire Invariant One Theorem \cite{kervaire}. In this paper, the authors needed to know that a particular Bousfield localization of equivariant spectra preserved commutative structure. My paper recovers and generalizes the theorem of Hill and Hopkins \cite{hill-hopkins} which provided conditions for such preservation to occur. Following Organizing Principle \ref{keyidea:build-model-categories}, my method of proving a preservation result is to try to put model structures on the category of objects with algebraic structure (e.g. commutative ring spectra). The following theorem reduces the question of preservation to a simpler question.

\begin{theorem}\label{bigthm}
Let $\M$ be a monoidal model category, $C$ a set of maps in $\M$, and $P$ an operad valued in $\M$. If $P$-algebras in $\M$ and in $L_C(\M)$ inherit (semi-)model structures such that the forgetful functors back to $\M$ and $L_C(\M)$ are right Quillen functors, then $L_C$ preserves $P$-algebras up to weak equivalence. For well-behaved $P$ there is a list of easy to check conditions on $\M$ and $C$ guaranteeing these hypotheses hold.
\end{theorem}

This theorem is proven in Section 3, and while the paper centers around this result (especially, checking its hypotheses) this is not where the hard work is being done. The proof really just involves a diagram chase and cofibrant/fibrant replacements, following Organizing Principle \ref{keyidea:replace}. A semi-model category is something slightly weaker than a model category, but which is much easier to build on a category of algebras and which still allows tools of homotopy theory (especially these (co)fibrant replacements) to be used. It satisfies all the axioms of a model 
category except the lifting of a trivial cofibration against a fibration and the factorization into a trivial cofibration followed by a fibration require the domain of the map in question to be cofibrant. The property of being cofibrant should be viewed as being analogous to being a CW complex or a projective module. Every object is weakly equivalent to a cofibrant object (following Organizing Principle \ref{keyidea:replace}) via the useful tool of cofibrant replacement. 
In order to obtain concrete, recognizable results, the paper specializes to two settings: where $P$ is a $\Sigma$-cofibrant 
operad and when $P$ is the Com operad. In the former case it has long been known how to transfer semi-model structures from $\M$ 
to $P$-algebras, and Theorem 5.1 recalls the procedure. In the latter case, my companion paper \cite{white-commutative-monoids} solved the problem and the main result of that paper is 
recalled in Theorem 6.2. In both cases the key point is that

\begin{keyidea} \label{keyidea1} In order to transfer a model structure to a category of algebras over a monad one must have good homotopical control over the free algebra functor. This often requires some kind of filtration so that free extensions (i.e. pushouts) in the category of algebras can be computed by some transfinite process in the underlying category $\M$.
\end{keyidea}

These filtrations often take many pages to develop, but they are not the main point of any such paper. They are more an artifact of the topologist's method of proof, which involves building a complicated machine to compute something via a transfinite process and then recovering a result of interest as a special case. The $\Sigma$-cofibrancy hypothesis effectively ensures good homotopical control in any cofibrantly generated model category $\M$. The Com operad is not $\Sigma$-cofibrant, but one can still obtain good control by making a hypothesis on $\M$. In \cite{white-commutative-monoids} the hypothesis is introduced as the Commutative Monoid Axiom, and has to do with the free commutative monoid functor $Sym(X)= S \coprod X \coprod X^{\otimes 2}/\Sigma_2 \coprod X^{\otimes 3}/\Sigma_3 \coprod \dots$ where $S$ is the monoidal unit and $\Sigma_n$ is the symmetric group on $n$ letters. The following is proven in \cite{white-commutative-monoids} and recalled as Theorem 6.2 in \cite{white-localization}:

\begin{theorem} \label{thm:commMonModel}
If a monoidal model category $\M$ satisfies the commutative monoid axiom (i.e. for any trivial cofibration $g$, the map $g^{\boxprod n}/\Sigma_n$ is a trivial cofibration) then commutative monoids inherit a semi-model structure from $\M$ which is a model structure if $\M$ satisfies the monoid axiom from \cite{SS00}.
\end{theorem}

The commutative monoid axiom is verified in \cite{white-commutative-monoids} for model categories of spaces, simplicial sets, chain complexes in characteristic zero, diagram categories, ideals of ring spectra, and positive variants of symmetric, orthogonal, equivariant, and motivic spectra. Theorems 5.1 and 6.2 cover the hypotheses of Theorem \ref{bigthm} about $P$-alg($\M)$ having a semi-model structure. In order to check the hypotheses about $P$-alg$(L_C(\M))$ we need in addition

\begin{keyidea} \label{keyidea:axioms-to-localization}
In order for Bousfield localization to preserve operad algebra structure one should verify that it respects the axioms of monoidal model categories, e.g. the Pushout Product Axiom, the Unit Axiom, the axiom that cofibrant objects are flat, the Commutative Monoid Axiom, and the Monoid Axiom.
\end{keyidea}

If $L_C(\M)$ satisfies all these axioms then Theorem \ref{bigthm} will prove that commutative monoids and algebras over $\Sigma$-cofibrant operads are preserved by the localization. Section 4 introduces a hypothesis on the maps to be inverted which guarantees the first three of these axioms hold for $L_C(\M)$. Indeed, Section 4 characterizes the localizations for which these axioms hold:

\begin{theorem} \label{thm:PP}
Assume $\M$ satisfies the pushout product axiom and that cofibrant objects are flat. $L_C(\M)$ satisfies these axioms (hence the unit axiom too) if and only if for all cofibrant $K$, all maps of the form $f\otimes id_K$ for $f\in C$ are weak equivalences in $L_C(\M)$. If $\M$ has generating cofibrations $I$ with cofibrant domains then it is sufficient to check the condition for $K$ in the set of (co)domains of $I$.
\end{theorem}

This theorem requires a fair bit of work to prove, but it is fun work for model category theorists. The case for cofibrant domains required the nifty Lemma 4.13 which I hope will help future users of model categories. With Theorem \ref{thm:PP}, Theorem 5.1, and Theorem \ref{bigthm} we have a list of checkable conditions so that localization preserves algebras over $\Sigma$-cofibrant operads. The conditions are checked for numerous examples in Section 5 and recover examples of Farjoun and Quillen for spaces and chain complexes respectively. Counterexamples are also given, including perhaps the first explicit example where the pushout product axiom fails to be satisfied for some $L_C(\M)$. Most of the work in the paper comes in checking the examples in Sections 5 and 7, since this is the only way I could tell if the hypotheses I introduced were good or not. This theorem demonstrates another key idea which is in the background of this work

\begin{keyidea}
The theory of monoidal categories can serve as a useful guide when proving results about monoidal model categories.
\end{keyidea}

In particular, a similar characterization to that in Theorem \ref{thm:PP} appeared in work of Brian Day. The condition precisely ensures that the localization respects the monoidal structure, and such localizations are dubbed Monoidal Bousfield Localizations in Section 4. Similarly, to check the commutative monoid axiom one must know that localization respects the functor Sym. This helps us check another condition in Key Idea \ref{keyidea:axioms-to-localization}, and appears in Section 6:

\begin{theorem}\label{thm:loc-preserves-cmon-axiom}
Assume $\M$ is a well-behaved monoidal model category satisfying the commutative monoid axiom. Suppose that $L_C(\M)$ is a monoidal Bousfield localization. Then $L_C(\M)$ satisfies the commutative monoid axiom if and only if $\Sym(f)$ is a $C$-local equivalence for all $f \in C$.
\end{theorem}

This result required a great deal of work to prove and was not satisfying to me because of the hypotheses required (currently hidden under the phrase ``well-behaved''). I hope to return to this result in the future and get a slicker proof without these hypotheses. With Theorems \ref{bigthm}, \ref{thm:commMonModel}, and \ref{thm:loc-preserves-cmon-axiom} we have achieved our goal of finding checkable conditions so that localization preserves commutative monoids. The conditions are checked for numerous examples in Section 7, including spaces, chain complexes, various models of spectra, and equivariant spectra. Results of Farjoun, Quillen, and Casacuberta, Gutierrez, Moerdijk, and Vogt are recovered and generalized, as well as new results for equivariant spectra.

Lastly, in order to make a complete story, Section 8 provides conditions on $\M$ and $C$ so that $L_C(\M)$ satisfies the monoid axiom. This is not necessary for Theorem \ref{bigthm} because semi-model structures suffice, but I felt any good theory of monoidal Bousfield localization should include results about the monoid axiom in case users of the paper need full model structures rather than semi-model structures. This section introduces a new axiom called $h$-monoidality, independently discovered in \cite{batanin-berger}, and checks it for a wide variety of model categories.

\section{Metaphors and Imagery}

Most of my work is phrased in the language of model categories, settings that allow abstract homotopy theory to be applied to numerous other fields. When I think of model categories, I don't think of the definition. Instead, I see a suite of commutative diagrams telling me exactly what I know I can do with a model structure. These diagrams seem beautiful and well structured. I have confidence that if I can just fit the pieces into the right places then they will commute (i.e. the various ways of moving through the diagram will agree) and prove what I need them to prove. Some of these diagrams are simplistic, e.g. triangles showing me I can factor maps. Others are complicated, e.g. towers to compute simplicial and cosimplicial resolutions (necessarily with many arrows between any two vertices in the tower), or large dimensional hypercubes which encode different possible orderings on letters in a word (these arise naturally when studying monoids). Simultaneously, I think of a suite of examples and counterexamples that warn me about which properties do not come for free, about times I've been surprised by the axioms (e.g. the model category of Graphs and its weak factorization systems), and about manifestations of model categories in radically different fields. 

\subsection{Examples of Model Categories}

The example I most often think about is that of topological spaces. When I inhabit this example I am immediately aware of safe and nice things I can do, such as CW approximation, mapping cylinders, finding liftings via sections to covering maps and properties of monomorphisms, and monoidal properties given by Cartesian product and the compact-open topology on hom spaces. I tend to picture spaces as things I can draw such as manifolds, but I also picture my favorite counterexamples as weakenings of the various manifold axioms. I like to think about the long line (which has trivial homotopy groups but is not contractible), the Sierpinski space (which demonstrates that Top cannot be locally presentable), the Sorgenfrey line, finite topological spaces, a pushout of an inclusion of compactly generated weak Hausdorff space which is not even injective, etc. My love of this subject began in point-set topology with Munkres's \textit{Topology} and for this reason I believe counterexamples are a fundamental and beautiful part of the subject. That said, I do at times zoom in to nice subcategories such as CW complexes, $\Delta$-generated spaces, and compactly generated spaces. When doing so I keep track of certain spaces I care about, mostly compact Lie groups and orthogonal groups $O(n)$.

My next most favorite examples come from algebra. I think of the projective model structure on Ch(R), and I have a mental switch which I can throw telling me whether we are thinking of bounded chain complexes (where everything is nice) or unbounded chain complexes where things can be tricky. When thinking about monoidal properties I have another switch which makes sure I am either working over a field of characteristic zero or proceeding with extreme caution. In this example I often switch from thinking about maps to thinking about objects: namely their fibers or cofibers. This shift makes things much easier, as I can simply ask if an object is acyclic rather than trying to use the functoriality of homology to study $H_*(f)$ for some map $f$. Similarly, the example of the stable module category is a great place to test conjectures about model categories because all objects are both fibrant and cofibrant, the category is stable so I can always reduce to studying objects (which are simply modules here), and yet I know many interesting examples that I can use to disprove conjectures. Example 4.1 is computed in this setting and as far as I know this is the first non-trivial example of a model category with a monoidal product which fails to satisfy the pushout product axiom or the monoid axiom.

After that I like to think of spectra. I envision them just like chain complexes, i.e. infinite chains with some strands connecting adjacent ones.  It is very easy to shift viewpoint to orthogonal spectra by simply allowing an action of $O(n)$ on the $n^{th}$ space $X_n$, and I view this action simply as an arrow from $X_n$ to itself.

\begin{figure}[ht!] \label{image:orthog-spectra}
\centering
\includegraphics[width=140mm]{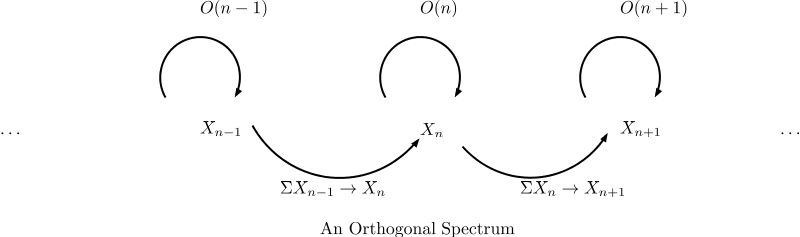}
\end{figure}

I have similar ways of thinking about equivariant spectra, now with actions of $G \times O(n)$, which I discuss in \ref{subsec:equivar}. Thinking of symmetric spectra is more difficult because they are sequence of simplicial sets rather than topological spaces. I have never felt as comfortable with simplicial sets as other model categories I study, partially because there are no good ways to visualize them (especially the face and degeneracy maps), so I always find myself thinking of simplicial complexes instead, i.e. of triangles, tetrahedra, and higher dimensional analogues. I know the various properties of simplicial sets and symmetric spectra, and can use them to make formal arguments. I rarely move inside to make point-set level arguments unless it is a lemma about smallness or presentability, which I understand by analogy to the category of sets. In section 7 of \cite{white-localization} I needed to make arguments at the level of simplices and it was difficult. In the end I relied on properties of the monoidal product, proving something holds for all simplicial sets $K$ rather than just for $K$ of the form $\partial \Delta[n]$ or $\Delta[n]$.

Other examples I often consider include simplicial presheaves, diagram categories with projective/injective/Reedy model structures, $W$-spaces, motivic symmetric spectra, and graphs. I think of the first three as categories of functors, effectively never using that $W$-spaces can be thought of as sequences because that would require fixing a skeleton for the category of finite CW complexes. I think of simplicial presheaves as diagrams where I know many formal properties are satisfied, but I never use simplex level arguments. I think of motivic symmetric spectra just like symmetric spectra, but at each level $n$ I see a shadow grading (so that the spectrum looks like a 2-dimensional lattice if viewed correctly). This tells me I can think of motivic symmetric spectra either as the stabilization of motivic spaces with respect to the functor $- \wedge \mathbb{P}^1$ or as the stabilization with respect to $-\wedge S^1_{alg}$ of the stabilization with respect to $-\wedge S^1_{top}$ of motivic spaces. Thinking by analogy to symmetric spectra lets me work in the motivic setting, but always cautiously to make sure my lack of background in algebraic geometry will not lead me to false conclusions. Lastly, when I think of the category of graphs many bits of mathematics I enjoy pop into my mind, e.g. Markov chains, electrical flows on graphs, algorithms to create spanning trees, etc. I need to shut these images down because to study the model category of graphs requires thinking about graph homomorphisms and zeta series on graphs, which I understand less about. The model category of graphs is surprising, and is best kept as an interesting example to study carefully on its own in some future work. For now it is too dimly lit and poorly explored to feel comfortable.

\subsection{Properties of Model Categories}

When I think generally, rather than in the context of specific examples, things become much easier for me. The model categorical tools I most often use feel comfortable and smooth. For instance, when I use cofibrant replacement I envision fattening up an object $X$ to its cofibrant replacement $QX$. I often view this as spreading $X$ out, e.g. shifting how I am viewing it so that instead of seeing a grid in 2 dimensions I see a 3 dimensional grid where the extra dimension was hidden from my previous perspective. I often think of $X$ as living inside of $QX$, though known counterexamples warn me not to use that intuition when writing down proofs. Instead, thinking this way makes me comfortable and confident, and this in turn increases the number of attempts I am willing to make on a problem in a single sitting without giving up. While fibrant replacement is formally dual to cofibrant replacement, it feels less natural to me. Instead of feeling smooth there seem to be jagged edges preventing me from feeling fully confident in its use. Instead of feeling white or silver like cofibrant replacement it feels darker and a bit shadowy. The main reason for this is that I always work in cofibrantly generated settings, and this means I have a collection of cells that let me build cofibrant replacement, e.g. CW approximation in the category Top of topological spaces, or projective resolution in a category Ch(R) of chain complexes. Injective resolution never felt as natural to me, and I learned about Top before simplicial sets (sSet) so I tend to like cofibrant replacement more. This preference has strongly affected my work, as in my thesis I focused on the ``cofibrant'' way to build model structures on algebras over a monad rather than the ``fibrant path object'' method. The benefit is that my work was able to apply to examples which had not previously been studied due to the fact that my theorems require different hypotheses (at times easier to satisfy) than theorems obtained via the fibrant approach.

When I think of properties a model category can have or fail to have, I again most often think of what these properties are good for rather than their literal meaning. For example, a model category $\M$ can be left proper, and this has a definition regarding the behavior of weak equivalences along pushouts by cofibrations. Rezk has given equivalent definitions in terms of behavior of functors on over/under categories. However, I think of left properness in terms of the two pictures below. The one on the left reminds me of Proposition 13.2.1 of \cite{hirschhorn} which says that in a left proper model category it's sufficient to test lifting against cofibrations between cofibrant objects. The one on the right reminds me of \cite{batanin-berger}'s proof that $\M$ is left proper if and only if all cofibrations are $h$-cofibrations. I have at times found this formulation easier to work with.

\[
\xymatrix{QX \ar[r] \ar@{^(->}[d] & X \ar[rr] \ar@{^(->}[d] & & A \ar@{->>}[d]^\simeq & & X \ar@{^(->}[d] \ar[r] \po & A \ar[r]^{\simeq} \ar[d] \po & B \ar[d] \\
QY \ar[r] \ar[urrr]^(.3){\exists} & Y \ar[rr] \ar@{..>}[urr]_{\therefore \exists} & & B & & Y \ar[r] & C \ar[r]_{\therefore \simeq} & D}
\]

When I learn a new fact about model categories I must find the right diagram to encode the fact; until I ``see" it presented in my own way I cannot believe the fact is true. Next, in order to remember the new fact I must fit it into my memory alongside all the other diagrams. I visualize this process as analogous to the way a computer writes to disk memory. My mind skims through all the facts I know and determines which are like this one. It then creates mental web strands to connect this new fact with the like facts identifies and fits the new diagram into its place at the barycenter of the related facts. I use this mental web frequently when searching my mind for workable proofs, and I share this imagery of ``latching new knowledge onto existing knowledge" with my students in every class I teach. I do a similar searching process when I prove something, to determine where precisely in the literature this new fact fits.

Let me give an example: when I started working with the pushout product axiom, that tells when a monoidal structure and a model structure are compatible, I frequently drew the picture below. It appears in the proof of Proposition 4.12, and is the key step to proving the main result of section 4, characterizing monoidal Bousfield localizations.

\begin{align*}
\xymatrix{K\otimes X \po \ar@{^{(}->}[r]^\simeq \ar[d] & K\otimes Y \ar[d] \ar@/^1pc/[ddr] & \\ L\otimes X \ar[r]^(.3){\simeq} \ar@/_1pc/[drr]_\simeq & (K\otimes Y) \coprod_{K\otimes X} (L\otimes X) \ar[dr]^{h\boxprod g} & \\ & & L\otimes Y}
\end{align*}

Another picture which often occurred was that of an $n$-dimensional punctured hypercube, e.g. the following for $n=2$ and $n=3$:

\[
\xymatrix{K \otimes K \ar[r] \ar[dd] & K\otimes L\ar[dd] &
K\otimes K \otimes K\ar[rr] \ar[dd] \ar[dr] & & K\otimes K \otimes L \ar[dr] \ar[dd] & \\
& & & K\otimes L\otimes K \ar[rr] \ar[dd] & & K\otimes L\otimes L \ar[dd] \\ 
L\otimes K \ar[r] & P & 
L\otimes K\otimes K \ar[rr] \ar[dr] & & L\otimes K \otimes L \ar[dr] & \\
& &  & L\otimes L\otimes K \ar[rr] & & P}
\]

In both cases $P$ is the pushout of the rest of the cube (called the `punctured cube') and it maps to $L^{\otimes n}:=L\otimes L \cdots L \otimes L$ because all vertices in the cube do so. When studying commutative monoids one must also take into account the $\Sigma_n$ action on the cube, which can equivalently be thought of as permuting the letters in the words which appear as vertices in the cube. These cubes appear all over Section 6, 7, and the Appendix to the companion paper \cite{white-commutative-monoids}. One difficulty of my mental imagery is that it's difficult to distinguish when a map (viewed as an edge in some diagram) respects the $\Sigma_n$ action and when it does not. This made working out the mathematics in Section 6 extremely difficult for me, and as a result I decided to simply convert the problem of finding a $\Sigma_n$-equivariant lift to the problem of finding any lift at all in a related diagram which I could understand better.

\subsection{Equivariant Spectra} \label{subsec:equivar}

In Sections 5 and 7 I apply the main results of the paper to several examples. I recover classical results about spaces, spectra, and chain complexes. The main new results are about $G$-equivariant spectra, where $G$ is a compact Lie group. In order to have a good monoidal product I work in the context of $G$-equivariant orthogonal spectra. I have a picture much like \ref{image:orthog-spectra} in mind, since to me a $G$-spectrum is just an orthogonal spectrum on which $G$ acts, i.e. an $\N$-graded sequence such that for all $n$, $G\times O(n)$ acts on $X_n$. I am aware that in every dimension $n$ there are also spaces $X_V$ related to $X_n$ but with a twisted $G$-action inherited from the $G$-action on the $n$-dimensional representation $V$. I view a $G$-spectrum as a chain with shadow versions of $X_n$ clustered in level $n$ for all $n$. I also see restriction and transfer maps between the shadow versions at different levels just like I see the suspension map data as floating ``between'' different levels. The structure maps are now a bit more complicated, because the maps $S^V \wedge X_W \to X_{V\otimes W}$ must be $G \ltimes (O(V)\times O(W))$-equivariant. However, this is easy to remember because $G$ acts by conjugation on maps, so I simply view it enveloping the orthogonal group actions.

In order to properly study commutativity for $G$-spectra, one must work with smash products indexed by $G$-sets. I often restrict attention to finite cyclic groups $G$ or symmetric groups so that I can see the action of $G$ on the smash product. I then write down proofs in the maximal generality possible and check that they hold for compact Lie groups $G$. In order to even define equivariant homotopy one needs the notion of $H$ fixed points for a subgroup $H$ of $G$. Here I again think of $G$ permuting points and it is easy in this light to see the points which are not permuted. I view subgroups via their permutation actions, and these are usually color-coded so that I can see which points they move and which they do not. At times I've had to work with the whole lattice of subgroups of $G$, or with families of subgroups of $G$ (i.e. sets of subgroups closed under conjugation and passage to subgroup). I again use color-coding for this (e.g. $H$ is usually red, $K$ is blue, and subgroups of $K$ are different shades of blue), and it makes it easy to keep track of the various family model structures (where you vary the weak equivalences according to which homotopy groups are seen to be isomorphisms by the family) by their colors. 

I picture the lattice of subgroups as another diagram, whose edges tell me whether one family is contained in another. When I draw this lattice I often draw it as a tower due to space constraints. One of the main results of this paper is that there are localizations which destroy some, but not all, of the commutative structure (equivalently, of the multiplicative norms). The structure which is not destroyed can be viewed as a change of family, so I often draw maps between these towers. This story is made rigorous using a collection of model structures on the category of $G$-operads (one for each feasible choice of families $\mathscr{F}_n$ of subgroups of $G\times \Sigma_n$) in ongoing joint work with Javier Guti\'{e}rrez.

Related to my picture of the lattice of family model structures is the lattice of universe model structures, where you vary the weak equivalences based on which $G$-representations are allowed (these model structures are discussed in \cite{hovey-white} among other places). I have often wondered about the relationship between these family model structures and these universe model structures. I envision a two dimensional grid of model structures on the category of ($\N$-indexed) $G$-spectra with some possible folding I have not yet understood which will tell me when a universe model structure is Quillen equivalent to a family model structure (i.e. encode the same homotopy theory). There is a natural way to relate universes and families via determining where isotropy subgroups embed, but I have not yet proven that this correspondence yields Quillen equivalences.

\subsection{Bousfield Localization}

When I think of localization, it is as a large arrow between categories satisfying a universal property. This arrow is a functor, and so specializes to small arrows for any objects in the category. Each such arrow takes an object to the closest object which is local, and by construction each such arrow will be a local equivalence. Thus, I think of the target of the arrow as a shade of the first, partially transparent, which is as much of the original as can be ``seen'' by the localization (i.e. by the maps which are being inverted). 

This paper is about monoidal localization, so I always assume the model category $\M$ has a monoidal product which respects the homotopy, i.e. descends to a monoidal product on the homotopy category $Ho(\M)$. This means that for any objects $X$ and $Y$, if I choose weakly equivalent objects $X'$ and $Y'$ from the respective fuzzy red clouds I see around $X$ and $Y$ then I will end up with a homotopy equivalence between $X\otimes Y$ and $X'\otimes Y'$. I want the same to be true after localization, i.e. in $L_C(\M)$ where $C$ is the set of maps I am inverting. This $L_C(\M)$ has more weak equivalences, so the clouds of equivalent objects have gotten bigger. I provide general conditions which guarantee that $L_C(\M)$ is again a monoidal model category, and then consider these conditions in my examples of interest. 

In spaces and bounded chain complexes the conditions are always satisfied. In spectra one needs to know that the localization is \textit{stable}, i.e. commutes with suspension. I visualize this as a map from a chain $X$ to its shade $L_C(X)$ (another chain) and I am simply asking that if I shift my window on $X$ to the right then shift it down to $L_C(X)$ then this will be the same as first going down to $L_C(X)$ then shifting to the right. I even hear the sound of gears locking together when I visualize these shifts. For equivariant spectra the conditions are more subtle. In the end it comes down to the functor Sym (which takes an object $X$ to the free commutative monoid on $X$) respecting local equivalences. If only some equivalences are respected (e.g. those in a family) then failure can occur. So I need my localization to respect all families, i.e. all colors, which is equivalent to asking it to respect the family of all subgroups. 

\subsection{Infinity Categories}

This paper does not deal with $\infty$-categories, but I often receive questions about them when I speak about this material and I want to share the imagery I associate with them. Currently, $\infty$-categories (in particular, quasicategories) are very popular among young researchers and some want to replace model categories with them. The upside would be that many statements of interest would be easier to prove, e.g. determining when a localization is monoidal, or when one has a good homotopy theory on a category of algebras. The downside is that one loses the ability to do many of the nice computational constructions one can with a model category, e.g. (co)fibrant replacement, (co)simplicial resolutions, computing homotopy (co)limits, finding explicit lifts, etc. People often ask why I don't work in the easier setting. The simplest answer is that in all my projects so far, $\infty$-categories were insufficient to encode the items I needed to study. In \cite{white-localization} the problem is that $\infty$-categories are not very convenient for discussing strict commutativity; they are better suited to $E_\infty$-algebra structure. As the examples in Section 5 and 7 show, this is often the wrong concept to study. In my papers with Michael Batanin, $\infty$-categories were insufficient to encode the notions of the project: that is, $n$-operads, cofibrant replacements of $n$-operads, and algebras in weak $n$-category (with no higher cells, rather than $(\infty,n)$-categories with a contractible choice of higher cells). In my work with Javier Guti\'{e}rrez, the concept of cofibrant replacement is crucial to constructing the operads we wish to study, and this concept does not exist for general $\infty$-categories. More generally than all of these examples, I have yet to find a problem that I cannot solve using model categories but could solve using $\infty$-categories.

I understand why some want to research $\infty$-categories, especially those with a strong background in simplicial sets. It's a context where you are virtually guaranteed of success as long as you can keep up with all that has been written by others on the subject. However, I prefer to continue to work in the setting of model categories, where things feel more hands-on to me, where I understand and use all the tools, and where I already have familiarity translating general results into specific examples. In addition, I believe working with model categories makes my work more applicable, since any $\infty$-category theorist can easily translate my results into their setting (whereas translating the other way is often difficult and not always possible) but for the majority, who understand model categories but not yet $\infty$-categories, my results are already in a language they can use. Lastly, I like to work in fields where there are actual counterexamples! This was part of what brought me into topology in the first place and discarding model categories would be letting go of some of my favorite counterexamples. 

A few weeks ago I was exploring an art museum and I didn't have enough time to see everything. I love art museums and have learned over the years that art which is carefully detailed appeals to me most strongly. In particular, landscapes and large scale paintings of human endeavor always catch my eye. If I am pressed for time I usually skip the impressionist gallery, but on this day the entire special exhibit was on impressionists. The museum had pieces from all the best masters, and I could finally understand why this form of art appealed to so many. The paintings I saw perfectly captured the light at various parts of the day, but completely blurred out the actual objects of study. I realized my taste for model categories and away from $\infty$-categories is the same as my taste in art. I like to see the details, the precise and intricate point-set level constructions, and the clever ways of fitting known diagrams together to give new facts. I am not satisfied by just ``the essence'' of the object or its characterization by a universal property, though I do see the artistic value in that approach.

Research in mathematics can be tumultuous. A mathematician is fundamentally a truth-seeker, but until the truth is found must hold conflicting possible truths in mind simultaneously (e.g. when deciding whether to look for a proof or a counterexample). This is especially true in graduate school, before one has properly developed the appropriate gut instinct for whether a statement is true or false, or for how difficult the statement will be to prove. Graduate school can also be difficult for other reasons: students don't know whether or not they'll be successful in completing the degree, whether or not they'll get a job afterwards, where they might end up living, or what other responsibilities life might send them. In addition, graduate school is a time to learn and obtain expert level knowledge in the discipline, but this can be difficult because in many cases there are too many references to possibly read and sometimes things are well-known but not written down anywhere.

In this section I will tell the story of my PhD thesis in the hope that it might help a future graduate student understand a bit of the process, especially the fact that every grad student faces similar challenges, doubts, and questions. This section also contains context for \cite{white-localization}, information about how some of the ideas arose, and tips for making the most out of your time in graduate school.

\section{Story of the Problem}

I began working on \cite{white-localization} in August 2011, and it was not my first PhD project. During my third year in graduate school (2010-2011) my advisor Mark Hovey gave me my first project, to extend his work with Keir Lockridge and compute the homological dimension of the real (resp. complex) K-theory spectrum KO (resp. KU). I spent months learning the requisite background but found myself completely unable to make headway on the problem. Hovey and Lockridge already knew in 2009 that $2\leq dim(KU) \leq 3$ and that $4\leq dim(KO)\leq 5$, and that remains the state of knowledge today (here dim means global dimension as a ring spectrum). 

As attempt after attempt failed I became frustrated and devoted increasing amounts of time to my master's research in computer science with Danny Krizanc, which was just getting underway. I found that research much easier, and for about a month I began all my meetings with Mark by confessing I had nothing new to report about K-theory but had managed to prove some fact or other about navigation algorithms for autonomous agents moving on a graph. It amazes me that Mark let me get away with this, but I am very grateful he did. The side-project in computer science helped me rebuild my confidence and without my computer science degree I would not have the job I have today. In July of 2011 we decided to make one more concerted effort to resolve the dimension problem, and I made the following conjecture

\begin{conjecture}
If a localization $R\to R[v^{-1}]$ is sufficiently nice then it cannot reduce global dimension by more than 1.
\end{conjecture}

This conjecture would resolve the question for both KO and KU at once, since the only way for them to differ by 1 would be if $dim(KU)=3$ and $dim(KO)=4$. This really felt like the right idea to me: it would be an elegant solution, it would involve formal arguments rather than spectral sequence computations, and it was borne out by examples in the setting of pure algebra (rather than the setting of ring spectra, sometimes referred to as ``Brave New Algebra"). Sadly, after more than a month of working on this as hard as I could, I had nothing to show and we decided to find another thesis project. It must have been this experience that made Mark realize that I should be working on a problem featuring localization.

Over the next month I began to work out various facts about Hovey's new theory of Smith ideals of ring spectra, which had only just been defined and therefore seemed a perfect thesis topic to me (almost nothing was known and there was no chance of competition). However, after attending a conference in Germany in August, 2011, Mark returned to Wesleyan, called me into his office and excitedly told me he had found my thesis problem. His exact words were: ``Something I thought always works turns out not to in an exciting new example. You are going to figure out why and to find conditions to make it work." The example is now Example 5.7 in \cite{white-localization} and Mark had learned it from a talk by Mike Hill. To summarize it:

\begin{example} \label{example:hill}
There is a stable localization of equivariant spectra which does not preserve commutative ring spectrum structure.
\end{example}

My paper \cite{white-localization} finds conditions on a model category and on a set of maps one wishes to invert (via Bousfield localization) so that this preservation does occur. Specializing to equivariant spectra, these conditions tell us which maps we can invert without losing commutativity, and demonstrate that the example is ``maximally bad'' in the sense that it destroys as much commutative structure as it possibly could while still being stable (here stable means with respect to the monoidal unit, not with respect to all representation spheres, and that's part of the problem). This was an excellent thesis problem: it had a concrete application at the end, related to one of the most exciting results in recent years (the Kervaire theorem of \cite{kervaire-arxiv}), and the solution allowed me to do some cool work in model categories of independent interest. For instance, it required me to work out when commutative monoids inherit a model structure (a problem that had been open for 15 years, since \cite{SS00}) and then to work out when Bousfield localization respects monoidal structure (a kind of join of Hirschhorn's book \cite{hirschhorn} on localization with the chapter of Hovey's book \cite{hovey-book} on monoidal model categories).

\subsection{Story of the Development}

In the fall of 2011 I simultaneously read \cite{hirschhorn}, figured out conditions so that localization would preserve the pushout product axiom, and (with help from Mark) figured out sufficient conditions so that the localization would preserve monoid structure (preservation for commutative monoid structure only came later). Over the winter I used a new condition Mark came up with regarding something he called \textit{homotopical cofibrations} to prove a result about preservation of the monoid axiom, though I felt the hypotheses were restrictively strong. I made pleasantly steady progress throughout the 2011-2012 year, with new results at almost every meeting with my advisor. Since he was Department Chair we met once every two weeks, sometimes with a formal meeting in between if he had time, or with hallway conversations about model categories whenever the opportunity arose.

During this year I organized the Wesleyan Topology Seminar and met many of the experts in the area. Everyone I invited responded favorably, and all who knew my advisor expressed their excitement at seeing him again. When they came I learned the pleasure of discussing research with experts from an array of backgrounds, and I also got to know my advisor better through dinners and conversations with the speakers; we began to become friends. Through these speakers I got new references to read, potential applications for my work, and their extremely valuable first impressions and flashes of insight on my research program. This experience led me to seek out experts whose papers I had read and engage them at conferences and through email. Those experiences in turn helped me hone my ``elevator pitch'' for my research: I learned the quickest way to describe my results, how to make the results sound interesting, which questions to expect, and which questions to ask so that I could move the research forward. 

In addition to conducting the research in \cite{white-localization} and organizing the seminar, during this time I also wrote my master's thesis in computer science, and began in the spring to give talks in various seminars (organizing a seminar is also a great way to get invited to speak in other seminars). I found myself somewhat exhausted going into our spring break. I sent my master's committee my 80 page thesis and then took a two week vacation in the south of France to visit my then girlfriend. We decided to take a long weekend in Barcelona in the middle of this vacation. 

I had recently read several papers by a well-known mathematician in Barcelona named Carles Casacuberta, who had done work related to localization and preservation of algebra structure years before. When I found this work in November of 2011, I worried that it might subsume my project, but Mark convinced me that they were different, for reasons which are now spelled out in \cite{white-localization} at several points. I wrote to Carles to ask if we could meet for coffee, figuring this would be another way to get feedback as with the Wesleyan Topology Seminar. He responded by inviting me to give a talk and then taking me out to a very fancy lunch with his postdoc (and former student) Javier Guti\'{e}rrez. The conversation went so well that Carles invited me to come back for the summer of 2013, and now both Carles and Javier are co-authors of mine (on different projects).

True to form, within a few minutes of the end of my talk, Carles's first instinct was spot on and gave me an idea which greatly influenced my research program. He remarked that my preservation result (Theorem 3.2 in \cite{white-localization}) was general enough to hold for colored operads, not just for associative and commutative monoids. I enthusiastically agreed even though at the time I had no idea what a colored operad was. It wasn't until a full year later that I really understood the story for operads (now Section 6.6 in \cite{white-thesis}) and another year after that till I understood the version for colored operads (which has appeared in \cite{white-yau}). 

In the spring of 2012 I focused on the situation for commutative monoids, and I discovered the commutative monoid axiom in May, just after defending my master's thesis. Marcy Robertson was our visitor that week and I distinctly remember Mark excitedly going down the hallway to check the new condition with her and see if it was likely to be satisfied in examples of interest. I spent the summer doing some extremely technical work related to proving that it suffices to check the commutative monoid axiom on generating (trivial) cofibrations (now appendix A of \cite{white-commutative-monoids}), finding out when Bousfield localization preserves the commutative monoid axiom (now Section 6 of \cite{white-localization}), and working out the generalization to operads following Carles's suggestion. This was the majority of the hard work in the thesis, and was extra frustrating because (at least for the localization result) the only proof I could work out included hypotheses I felt sure were not necessary. 

In addition, it turned out to be subtle to find the right condition so that categories of operad algebras inherit a model structure. Mark has no papers featuring operads, so we had to learn the field together. A visit by John Harper in the middle of the summer helped convince me I had the right condition, but it was difficult to write down a human-readable proof. Mark insisted on having such a proof, and I learned a great deal about writing as he rejected three versions before I produced one which he was happy with. This last proof and explanation is how I have presented the result ever since: if $P$ is an operad in a monoidal model category $\M$ then in order to know that $P$-algebras inherit a model structure a cofibrancy price must be paid on either $P$ or on $\M$. The most general form of this is now in \cite{white-yau} and it recovers all results of this sort (i.e. about inherited model structures on $P$-algebras) while also proving new ones about operads which are levelwise cofibrant, a situation where the cofibrancy price is paid partially by $P$ and partially by $\M$.

\subsection{The writing process}

I planned in the fall of 2012 to apply for postdoctoral positions, write up the results from 2011-2012, give talks in various seminars, and work out examples of the theory I had developed, especially to the case of equivariant spectra so that I could understand Example \ref{example:hill} above. However, September of 2012 turned out to be one of the worst months of my life personally. My father was diagnosed with a dangerous form of cancer, my family's financial situation degraded rapidly for a different reason, an old injury in my shoulder returned and ended my ability to play volleyball (till then a passion equal to mathematics in my life), and a long term relationship ended. Suddenly I found myself needing to return to Chicago frequently, spending huge amounts of time dealing with financial matters, and doing my best to provide support for my family despite having no solid ground on which to stand.

Mark was extremely supportive during this time. I lost my ability to focus on mathematics, and our regular meetings often turned to personal matters. He helped me realize it is okay to have periods like this from time to time; perhaps they should even be expected. He also helped me find ways to get to Chicago and I will never stop being grateful for his support during this time. The only research I accomplished in the fall was extending my results from the context of model categories to semi-model categories. At the time this seemed trivial to both me and him, but it turned out to be important in my development as a mathematician. Semi-model categories have appeared in the majority of my papers, often in places where model structures do not exist. Although Mark invented them (in \cite{hovey-monoidal}), he does not trust them much and made me carefully show him all steps of every result I claimed. In hindsight, it makes sense that I developed this extension at this time, since I had just read Markus Spitzweck's thesis \cite{spitzweck-thesis} the previous summer where semi-model categories were first explored in depth. From the point of view of this User Guide, semi-model categories are the reason Section 8 of \cite{white-localization} is not central to the story, and the reader can still have preservation results without needing to digest the meaning of \textit{homotopical cofibration} (now called $h$-cofibration following \cite{batanin-berger} who independently discovered several results in Section 8 in 2013).

I spent winter break supporting my family, and by the start of the spring semester things had stabilized. I had a massive backlog of writing to do, and still needed to work out the application to equivariant spectra. I spent the spring writing the parts of the story I understood best (from 2011 mostly) and trying to learn enough about equivariant spectra to finish my thesis. I got many helpful references from Carolyn Yarnall and Kristen Mazur at a conference that spring, but the word on the street was that the definitive reference would be the appendix to \cite{kervaire-arxiv} which had not appeared yet. When I mentioned this to Mark he shared his own view of equivariant spectra via model categories and using that foundation I was able to work out the fact that all my model category axioms were satisfied by the positive stable model structure on equivariant orthogonal spectra. When I asked him the next fall for a reference I could cite we ended up writing \cite{hovey-white} together to fill this gap in the literature, though that work has now been subsumed by \cite{kervaire-arxiv}.

The appearance of \cite{hill-hopkins} provided the last piece of the puzzle, as there were now numerous equivalent conditions a set of morphisms could satisfy, and these conditions implied preservation of commutative structure for $G$-spectra, at least for $G$ a cyclic $p$-group and for a particular kind of localization. I proved that the conditions in my Theorem 6.5 implied one of the conditions in \cite{hill-hopkins} and hence that my more general setting ($G$ could be a compact Lie group, and I was inverting a set of maps instead of a single homotopy element) included as a special case the result required for the Kervaire paper, bringing my thesis problem full circle.  From here it was also easy to see why Example \ref{example:hill} failed to satisfy these hypotheses, i.e. it is not a counterexample to my main theorem. I met Mike Hill at a conference in April and he kindly checked step by step the application of my general theorem to equivariant spectra. When he told me it all looked in order I knew my thesis was finished.

At this same conference I began a new project with Aaron Mazel-Gee and Markus Spitzweck in motivic homotopy theory. While this project never came to fruition, it was the starting point for the work I did that summer with Casacuberta in Barcelona. From April of 2013 till September of 2013 I tried to balance writing up the results in \cite{white-localization} with taking on new projects, giving talks in seminars, and preparing for the job market. I learned that it's always more fun to explore new math than to write up things I already understand, and for this reason I have to be very disciplined to actually write things up.

I spent May and June visiting Carles, and we proved the existence of a new localization in the context of motivic symmetric spectra, then began working out applications of this result. I also began a project on equivariant operads with Javier which we finally finished in the summer of 2015. In June I traveled to Nice (France), Nijmegen (Holland), and Lausanne (Switzerland) to give talks and discuss mathematics. My visit to Clemens Berger in Nice led to me strengthening \cite{white-localization} in two ways. First, he correctly suggested that the main result in Section 4 should be an ``if and only if'' statement. Second, he showed me a few things about $h$-cofibrations I didn't know and this led to cleaner statements and proofs in Section 8. During this visit I also had my first insight into the question of lifting left Bousfield localizations to categories of algebras, i.e. applying localization to commutative monoids (studying $L(CAlg(\M))$) rather than only looking at commutative monoids after localizing (studying $CAlg(L(\M))$). In the fall of 2013 I wrote a grant application to explore this problem with Michael Batanin and it ended up funding my trip to Australia in the summer of 2014.

I returned to the US in July to teach a short course in statistics, then returned to Barcelona in September using NSF funds for a conference on homotopy type theory. There Carles and I finished the first draft of our paper, though it took a long time until we could work that into a version to submit. All in all, the summer of 2013 was an extremely productive time for me and I found it very intellectually stimulating. I am grateful to Carles for hosting me and for making it possible. Once I returned from Barcelona in September, I spent October-January applying for jobs, doing interviews, giving talks, and writing \cite{hovey-white} with Mark so that I had a framework to write \cite{white-localization}. 

I accepted my job with Denison in early February and then turned my various typed-up results, lecture notes from my talks, and hand-written notes from meetings with various mathematicians into \cite{white-commutative-monoids} and \cite{white-localization} in March of 2014. I spent April writing my thesis and then spent May-August in Australia working with Michael Batanin. This was another intellectually stimulating time and we planted the seeds for at least three papers. Unfortunately, my Australia trip ended the same day that my job at Denison began and I learned first-hand that the first semester teaching takes up all of a new faculty member's time. In hindsight, I should have submitted my papers in March or April, because the publication process takes a long time and the delay in submitting did not improve the papers at all. My advice for graduate students is to send the paper to interested experts, spend two weeks making improvements based on their feedback (if any), put it on arXiv, and then submit it two weeks later unless you get more feedback. Although it seems like putting things on arXiv and submitting them are big commitments, it's always possible to make changes, minor changes are basically expected, and referees will make you change things anyway. So there is just no good reason to delay the process.

\section{Abstract for general audience}

Mathematics at its most basic level is the study of abstract thinking. Category theory follows this approach, and interprets all branches of mathematics as the study of objects and the relationships between them. For example, the objects might be people and you might say two are related if they're friends, or if one follows the others Twitter feed, or if they were in the same graduating class of high school, etc. Or the objects might be real numbers and you might decide one is related to another if it's smaller. Or the objects might be shapes and you could say two are related if they have the same symmetries (so a pentagon and a 5-pointed star are related, but the pentagon and the square are not related). When two objects are so tightly related that we want to view them as the same in all settings then we say those objects are equivalent (mathematicians use the term \textit{isomorphic}). For example, we might reasonably decide that all squares of side length 2 inches are the same. On the other hand, maybe someone else wants to study squares that come equipped with a color, so that red squares can be distinguished from blue squares. That's a different category, in which two squares with the same color and same side length would be considered equivalent. All of the work in this paper is in the setting of category theory, but to understand it we must now introduce another player.

Localization is a fundamental tool in mathematics that allows one to zoom in on the pertinent information in a problem. In the context of category theory, localization is a way to view two different (i.e. non-equivalent) objects as equivalent, e.g. deciding we don't care any more about color and now a red square and a blue square can be equivalent if they have the same side length. What's happening here is that we are putting on different eyeglasses when looking at the objects we want to study. Mathematically it means we're allowing more relationships between the objects, e.g. allowing blue and red to be related when before they were not. As humans we do this all the time. For example, if two driving routes take the same amount of time we might view them as equivalent. If two types of pasta in the supermarket cost the same we might view them as the same (if we don't care that much about pasta), and that's a valuable way to focus in on the information (in this case price) that really matters to us. 

In order to best study this localization procedure, we work in the setting of special kinds of categories called model categories. These are categories that come equipped with a specific localization we plan to do but have not done yet, i.e. a specific collection of relationships we want to eventually view as identifying the objects they are relating (these relationships are called \textit{weak equivalences} because they are not equivalences yet). For example, our category might be the category whose objects are shapes of all colors and where one shape is related to a second shape if they have the same color and if all the symmetries of the first are also symmetries of the second (e.g. a triangle is related to a hexagon). Then two objects are going to turn out to be equivalent if they have the same shape, size, and color. The weak equivalences could be relationships that ignore side length, so that two objects are weakly equivalent if they have the same shape and color, but not necessarily the same size. 

Model categories admit a special kind of localization called Bousfield localization (named after the mathematician Pete Bousfield) that transforms a model category into another model category with even more weak equivalences (i.e. we plan to view even more objects as equivalent). This procedure sends every object X to a closely related one (we'll denote it L(X)) that is equivalent to the original object according to the new notion of equivalence, but not according to the old notion. In our example above of shapes, sizes, and colors the new weak equivalences could be maps where the symmetries of one shape are symmetries of the second, but with no mention of side length or color (so we've added more relationships). The result is that two objects are weak equivalent according to the new weak equivalences if they have the same shape, but not necessarily the same size or color.

This paper is fundamentally about Bousfield localization. I studied how much structure on an object X is destroyed by the passage to L(X). Specifically, I was interested in algebraic structure on X. To a mathematician, algebra is a powerful computational tool and a great way of determining whether two objects are equivalent or not. Algebraic structure on an object should be thought of like icing on a cake. If two cakes have different icing then that's one sure-fire way to know they are different. Going back to our example of shapes, the information regarding the symmetries of the shape can be viewed as algebraic structure. One can describe a shape by the number of sides it has and where it's located in space, but this information about symmetries is extra and is often very useful. It's one way we can tell two triangles apart, for example (e.g. if one is equilateral and the other is obtuse). 

One way to encode algebraic structure in a category (i.e. to allow objects to possess algebraic structure) is via gizmos called \textit{operads}. For a given type of algebraic structure you want to study (e.g. a flavor of icing on the cake) these gizmos tell you exactly which objects have that structure. If you've got an object X that has the structure it's a natural question to ask whether its localization L(X) still has that structure. This paper answers that question in general by writing down exactly what must be satisfied in order for the algebraic structure to be preserved. It then goes on to work out specifically what that answer means in a number of model categories people have studied. I want to pause for a moment to explain why it was important to work out these specific cases of the general result.

Perhaps unsurprisingly, not all mathematicians like category theory. I love it, because if you can prove a theorem in category theory then it's true in all branches of mathematics. Similarly, if you prove a result about model categories then it's true in every specific model category out there. However, if someone is working in a specific branch of math (e.g. geometry) and the theorem they were trying to prove is proven using category theory they might be understandably frustrated; it feels like cheating. The relationship between category theory and the rest of math is much like the relationship between math and the rest of science. That is to say: category theorists produce results which can help in all areas of math, if someone takes the time to translate from the category theoretic jargon into the language folks in those areas are used to. Similarly, mathematics produces results useful all over science: both in physical science and social science, but often researchers in those fields don't want to go and learn a bunch of math in order to understand what has been done, so it falls to interdisciplinary researchers (e.g. mathematical physicists, mathematical biologists, mathematical economists, etc.) to bridge the gap and translate these general results into specific results in those fields. It often happens in economics that a team consisting of a mathematician and an economist are jointly awarded a Nobel Prize, because the mathematician worked out the general theory while the economist found lots of real-world applications for it. I think the analogy can be taken one step further. The relationship between a mathematician and a general scientist is like the relationship between the person who writes a cookbook and the actual cook in the kitchen. Alone neither might be successful but together they can produce something to better the world. If cookbook authors went around using all sorts of jargon that chefs could not recognize (e.g. discussing the use of ``positively curved metallic tools controlled via torque" instead of calling them ``spoons") they'd be doing their own work a disservice.

In my case I wanted people to use my work. So I learned the specific properties of a number of different model categories, and in each case I made sure they hypotheses for my general result held in those settings. This way I knew my work applied to a number of subfields of math, including algebra (the study of algebraic structure), topology (the study of space), stable homotopy theory (the study of when two spaces can be continuously deformed to become equivalent), representation theory (a way of studying objects based on how they cause other objects to change), and even to category theory itself. In each of these settings I learned the domain-specific jargon and proceeded to state my main result in that language, in the hopes that researchers in those fields would be able to use my work and comfortable doing so.



\begin{thebibliography}{1}

\bibitem[BB13]{batanin-berger}
Michael Batanin and Clemens Berger.
\newblock Homotopy theory for algebras over polynomial monads, preprint
  available electronically from http://arxiv.org/abs/1305.0086, to appear in \textit{Theory and Application of Categories}.

\bibitem[Bou75]{bousfield-localization-spaces-wrt-homology}
A.~K. Bousfield.
\newblock The localization of spaces with respect to homology.
\newblock {\em Topology}, 14:133--150, 1975.

\bibitem[Bou79]{bous79}
A.~K. Bousfield.
\newblock The localization of spectra with respect to homology.
\newblock {\em Topology}, 18(4):257--281, 1979.

\bibitem[Day73]{day-localisation} Brian Day. 
\newblock Note on monoidal localisation. 
\newblock Bull. Austral. Math. Soc. 8 (1973), 1--16


\bibitem[GZ67]{gabriel-zisman} Peter Gabriel, Michel Zisman. 
\newblock Calculus of fractions and homotopy theory, 
\newblock Springer, New York, 1967. Ergebnisse der Mathematik und ihrer Grenzgebiete, Band 35.

\bibitem[HH13]{hill-hopkins}
Michael~A. Hill and Michael~J. Hopkins.
\newblock Equivariant multiplicative closure, preprint available electronically
  from http://arxiv.org/abs/1303.4479.
\newblock 2013.

\bibitem[HHR11]{kervaire}
Michael~A. Hill, Michael~J. Hopkins, and Douglas~C. Ravenel.
\newblock A solution to the {A}rf-{K}ervaire invariant problem.
\newblock In {\em Proceedings of the {G}\"okova {G}eometry-{T}opology
  {C}onference 2010}, pages 21--63. Int. Press, Somerville, MA, 2011.

\bibitem[HHR15]{kervaire-arxiv}
Michael~A. Hill, Michael~J. Hopkins, and Douglas~C. Ravenel.
\newblock On the non-existence of elements of kervaire invariant one, version 4, preprint available electronically from http://arxiv.org/abs/0908.3724. To appear in {\em Annals of Mathematics}.

\bibitem[Hir03]{hirschhorn}
Philip~S. Hirschhorn.
\newblock {\em Model categories and their localizations}, volume~99 of {\em
  Mathematical Surveys and Monographs}.
\newblock American Mathematical Society, Providence, RI, 2003.

\bibitem[Hov99]{hovey-book}
Mark Hovey.
\newblock {\em Model categories}, volume~63 of {\em Mathematical Surveys and
  Monographs}.
\newblock American Mathematical Society, Providence, RI, 1999.

\bibitem[Hov98]{hovey-monoidal}
Mark Hovey.
\newblock Monoidal model categories, preprint available electronically from
  http://arxiv.org/abs/math/9803002.
\newblock 1998.

\bibitem[HL10]{hovey-lockridge} Mark Hovey and Keir Lockridge.
\newblock Homological Dimensions of Ring Spectra. {\em Homology, Homotopy, and Applications} 15(2), 2013, pages 55-71.

\bibitem[HW13]{hovey-white}
Mark Hovey and David White.
\newblock An alternative approach to equivariant stable homotopy theory,
  preprint available electronically from http://arxiv.org/abs/1312.3846.
\newblock 2013.

\bibitem[MMWWY15]{user-guide-project} Cary Malkiewich, Mona Merling, David White, Frank Lucas Wolcott, Carolyn Yarnall, The User's Guide Project: Giving Experiential Context to Research Papers, {\em Journal of Humanistic Mathematics}, Vol. 5, Issue 2 (July 2015), pages 186-188. .

\bibitem[SS00]{SS00}
Stefan Schwede and Brooke~E. Shipley.
\newblock Algebras and modules in monoidal model categories.
\newblock {\em Proc. London Math. Soc. (3)}, 80(2):491--511, 2000.

\bibitem[Spi01]{spitzweck-thesis}
M. Spitzweck,  Operads, algebras and modules in general model categories, preprint available electronically from http://arxiv.org/abs/math/0101102. 2001.

\bibitem[Whi14a]{white-commutative-monoids}
David White.
\newblock Model structures on commutative monoids in general model categories,
  available as arxiv:1403.6759. {\em J. Pure Appl. Algebra}, 221 (2017), no. 12, 3124-3168.

\bibitem[Whi14b]{white-localization}
David White.
\newblock Monoidal bousfield localizations and algebras over operads, available
  as arxiv:1404.5197.
\newblock 2014.

\bibitem[Whi14c]{white-thesis}
D. White, Monoidal bousfield localizations and algebras over operads. 2014. Thesis (Ph.D.)-Wesleyan University.

\bibitem[WY15]{white-yau} David White and Donald Yau.
\newblock Bousfield Localizations and Algebras over Colored Operads, available as arXiv:1503.06720. To appear, {\em Applied Categorical Structures}.

\end{thebibliography}

\end{document}